\documentclass[12pt]{article}

\usepackage[latin1]{inputenc}
\usepackage[english]{babel}
\usepackage{amsmath}
\usepackage{amssymb}
\usepackage{amsxtra}
\usepackage{amstext}
\usepackage{amsthm}
\usepackage{graphicx}

\theoremstyle{definition}
  \newtheorem{defin}{Definition}[section]
  \newtheorem{nota}{Remark}
  
\theoremstyle{plain}
  \newtheorem{teor}{Theorem}[section]

  \newtheorem{lema}{Lema}[section]

\newcommand{\re}{\mathbb{R}}

\title{On the Riemann problem for the simplified Bouchut-Boyaval system}
\author{Richard De la cruz\thanks{On license from Universidad Pedagógica y Tecnológica de Colombia - UPTC. \emph{e-mail:} {\tt richard.delacruz@uptc.edu.co} or {\tt radelacruzg@unal.edu.co}} \\ Juan Galvis\thanks{\emph{e-mail:} {\tt jcgalvisa@unal.edu.co}} \\ Juan Carlos Juajibioy\thanks{\emph{e-mail:} {\tt jcjuajibioyo@unal.edu.co}} \\ Leonardo Rendon\thanks{\emph{e-mail:} {\tt lrendona@unal.edu.co}}}
\date{Departamento de Matemáticas, Universidad Nacional de Colombia.}

\begin{document}


 \maketitle

 \begin{abstract}
In this paper we study the one-dimensional Riemann problem for 
a new hyperbolic system of three conservation laws 
of Temple class. This systems  it is a simplification of a recently propose  system of five conservations 
laws by Bouchut and Boyaval that model viscoelastic fluids. 
An important issues is that the  considered $3\times 3$ system is such that 
every characteristic field is 
linearly degenerate. 
Then, in despite of the fact that it is of Temple class,   the analysis of the Cauchy problem is more involved since general results for such 
a systems are not yet available. We show a explicit solution for the Cauchy problem with
initial data in $L^\infty$. We also study the Riemann problem for this system.
Under suitable generalized Rankine-Hugoniot relation and entropy condition, both existence and
uniqueness of particular delta-shock type solutions are established.
\end{abstract}
\section{Introduction}
The modeling of viscoelastic materials and fluids it is important for many applications. In particular,
a viscoelastic fluid is a material that exhibit both viscous and elastic characteristics upon deformation.
Examples of viscoelastic fluids that are important for applications are: latex paint,  gelatin, unset cement, liquid acrylic, asphalt and biological fluids such as synovial fluids, among others.
In \cite{Bouchut} the authors introduced a  new system of 
conservation laws that  models shallow viscoelastic fluids. This new system is motivated by Bouchut and Boyaval in
\cite[eq(5.6)]{Bouchut} and is written as 
\begin{equation*} 
 \begin{cases}
  \rho_t + (\rho u)_x=0, \\
  (\rho u)_t+ (\rho u^2+\pi)_x=0, \\
  (\rho \frac{\pi}{s^2})_t+ (\rho u \frac{\pi}{s^2}+u)_x=0, \\
  s_t+us_x=0, \\
  c_t+uc_x=0,
 \end{cases}
\end{equation*}
where $\rho$ denotes the layer depth of fluid, $u$ is the horizontal velocity, $s$
is related to the stress tensor and it is a conserved quantity,
$\pi$ is the relaxed pressure and $c > 0$ is introduced in order to parametrize the speeds.
This system describes a simple model for a thin layer of non-Newtonian viscoelastic fluid over a given topography at the bottom when the movement is driven by gravitational forces such as geophysical flows (mud flows, landslides, debris avalanches).

In \cite{Lu}, since $s$ is a conserved quantity, the author  considers the case $s=const. >0$.  Additionally,
we observe that the field $c$ does not appear in  the first four equations and, 
in order to simply even further, introduce the new variable $v=\frac{\pi}{s^2}$. 
After this observations, it is obtained the following simplified viscoelastic shallow fluid model,
\begin{equation} \label{sistema}
 \begin{cases}
  \rho_t + (\rho u)_x=0, \\
  (\rho u)_t + (\rho u^2+s^2 v)_x=0, \\
  (\rho v)_t + (\rho u v + u)_x=0.
 \end{cases}
\end{equation}
We refer to the system above as the \emph{simplified Bouchut-Boyaval system} (sBB) and it was introduced 
by Lu in \cite{Lu} as simplified version 
of a model proposed by Bouchut and Boyaval. 
We note that this  system is of Temple class and therefore, it is also it is of Rich type. 
We note that sBB is of Rich type but it is not diagonal, so its analysis is not standard. 
In this paper we are concern with the Riemann and Cauchy problems for (\ref{sistema}).
The existence of global weak solutions including vacuum regions, was obtained in  \cite{Lu} using the vanishing viscosity method in conjunction  with the compensated compactness argument. We also mention that, when dealing with the system \eqref{sistema},   one of the main difficulties is to obtain  existence and uniqueness of solutions of Cauchy problems in  the  presence of vacuum regions, that is, regions where the layer deep $\rho =0$.

We also note that, there are numerous studies on existence and uniqueness for general Temple class system \cite{Baiti-Bressan,Bianchini1,Bianchini2,Bressan-Goatin,heibig,Serre2}. 
However, some of these results do not apply to \eqref{sistema} since
it has all fields being
linearly degenerate and the initial data may have oscillations.\\
In this paper we obtain explicit solutions for a Cauchy problem associated to the sBB system \eqref{sistema} 
with $\rho_0(x)\ge \underline{\rho}>0$
with a possibly oscillating initial data. Also, we construct the Riemann solution for the system focussing our attention on delta shock
waves of certain type. The existence and uniqueness of solutions involving delta shock waves can
be obtained by solving the generalized Rankine-Hugoniot relation under a entropy condition \cite{Danilov-Shelkovich,Li-Zhang}.\\
The paper is organized as follows, in Section 2 we present the problem and put conditions on the initial 
data for physical properties are maintained in $\rho$.
In Section 3, we show the explicit solutions for the Cauchy problem associated with the sBB system \eqref{sistema} without the presence of vacuum regions. In Section 4, we solve the Riemann problem and we observe that the first and third contact discontinuity are asymptotic to the vacuum.
In the last section, we study the existence and uniqueness of solutions delta shock waves type.

\section{Properties of the  simplified  Bouchut-Boyaval system and some assumptions}
The eigenvalues associated to the system \eqref{sistema}
are given by,
\begin{equation} \label{VPropios}
 \lambda_1 = u-\frac{s}{\rho}, \quad \quad \lambda_2=u \quad \mbox{ and } \quad \lambda_3 = u+\frac{s}{\rho},
\end{equation}
where the corresponding Riemann invariants are 
\begin{equation}\label{invariantesRiemann}
 R_1=s^2 v-su, \quad \quad \quad R_2=v+\frac{1}{\rho} \quad \mbox{ and } \quad R_3=s^2 v+su.
\end{equation}
From  the expressions for the eigenvalues and the Riemann invariants we obtain
$$ \lambda_1 = \frac{R_3}{s}-sR_2, \quad \quad \lambda_2=\frac{1}{2s}(R_3-R_1) \quad \mbox{ and } \quad \lambda_3=sR_2-\frac{R_1}{s}. $$
From here we can see that system \eqref{sistema} is linearly degenerate.
On the other hand,   we have that for each $i, j,k \in \{ 1,2,3 \}$ with $j \neq i$, $k \neq i$, it holds
\begin{equation}
 \frac{\partial}{\partial R_j} \left( \frac{\frac{\partial \lambda_i}{\partial R_k}}{\lambda_k-\lambda_i} \right) = 
 \frac{\partial}{\partial R_k} \left( \frac{\frac{\partial \lambda_i}{\partial R_j}}{\lambda_j-\lambda_i} \right).
\end{equation}
This means that system \eqref{sistema} is of \emph{Rich} type. We recall that this classifications if due to \cite{Serre1}.\\

In this manuscript we focus on the study of the LBB system of conservation laws \eqref{sistema} with bounded initial data
\begin{equation} \label{datoincial}
\begin{aligned}
 & (\rho(x,0),u(x,0),v(x,0))=(\rho_0(x),u_0(x),v_0(x)), \quad x\in\re  \\ & \rho_0(x) \ge \underline{\rho}=const. > 0,
\end{aligned}
\end{equation}
subject to the following  conditions:
\begin{enumerate}
 \item[H1:] The functions $\rho_0$, $u_0$ and $v_0$ satisfy
            \begin{equation} 
             \begin{aligned}
              c_1 \le u_0(x)-sv_0(x) \le c_2, \quad c_3 \le u_0(x)+sv_0(x) \le c_4, \\
              \text{and }v_0(x)+\frac{1}{\rho_0(x)} > c_5,
             \end{aligned}
            \end{equation}
where $c_i$, $i=1,\dots,5$, are suitable constants satisfying
\begin{equation}
 c_5-\frac{c_4-c_1}{2s} > 0.
\end{equation}
 \item[H2:] The total variations of $u_0(x)-sv_0(x)$ and $u_0(x)+sv_0(x)$ are bounded.
\end{enumerate}

The conditions H1 and H2 are somehow natural to impose since they ensure that $\rho$ is positive giving a physical meaning to the sBB system \eqref{sistema}.
 
As we mentioned before, we note that in \cite{Lu}, the author shows existence of solutions for the Cauchy problem \eqref{sistema}-\eqref{datoincial} for the case  $\rho_0(x) \ge 0$.  This is done 
using the vanishing viscosity method  and a compensated complicatedness  argument.
In \cite{Lu} it is also shown that all entropies associated to \eqref{sistema} are of the form,
\begin{equation} \label{entropia}
 \eta(\rho,u,v)= \rho \left( F(u+sv) + G(u-sv)+H(v+\frac{1}{\rho}) \right),
\end{equation}
where $F,G,H$ are arbitrary functions having entropy flux
\begin{equation}   \label{flujoentropia}
 q(\rho,u,v)=(\rho u+s)	F(u+sv)+(\rho u-s)G(u-sv)+\rho u H(v+\frac{1}{\rho}).
\end{equation}
Moreover, if the functions $F,G$ y $H$ are convex, then, the entropy is also convex (see  \cite[Theorem 2]{Lu}). 
Thus, from each convex pair $(\eta,q)$  we have the following condition 
\begin{equation} \label{condEntropia}
 \eta_t(\rho,u,v)+q_x(\rho,u,v)=0
\end{equation}
in the sense of distributions.


\section{Explicit solutions}

In  this section we obtain explicit solutions for the Cauchy problem associated to sBB system with initial data \eqref{datoincial} subject to H1-H2 conditions. 
For this purpose, we used the results given Wagner, Weinan and Kohn, Li, Peng and Ruiz (for example see \cite{Wagner,Weinan-Kohn,Li-Peng-Ruiz,Peng2}).

 We uses the Euler-Lagrange (E-L) transformation $(t,x) \to (t,y)=(t,Y(t,x))$ defined by
$$ dy= \rho \, dx - \rho u \, dt \quad \text{ and } \quad Y(0,x)=Y_0(x) \stackrel{def}{=} \int_0^x \rho_0(\xi) \, d\xi . $$
In Lagrangian coordinates, the system \eqref{sistema} becomes
\begin{equation}\label{sistemaLin}
 \begin{cases}
  \omega_t-\nu_y=0, \\
  \nu_t+s^2\kappa_y=0, \\
  \kappa_t+\nu_y=0,
 \end{cases}
\end{equation}
where  $\omega$ denotes the quantity  $\frac{1}{\rho}$ in Lagrangian coordinates, that is, 
$\omega(t,y)=\frac{1}{\rho(t,x)}$, and we also have $\nu(t,y)=u(t,x)$ and $\kappa(t,y)=v(t,x)$.\\
The eigenvalues associated to  \eqref{sistemaLin} are given by
\begin{equation} \label{valoresPropiosLin}
 \widetilde{\lambda}_1=-s, \quad \widetilde{\lambda}_2=0, \quad \widetilde{\lambda}_3=s,
\end{equation}
and the the corresponding Riemann invariants are given by
\begin{equation}
 R_1=s^2 \kappa -s \nu, \quad R_2=\nu+\omega, \quad \text{ and } \quad R_3=s^2 \kappa -s \nu.
\end{equation}

In Lagrangian coordinates the entropy condition \eqref{condEntropia} transforms into 
\begin{equation}
 \widetilde{\eta}_t(\omega,\nu,\kappa) + \widetilde{q}_x(\omega,\nu,\kappa)=0,
\end{equation}
for each  $\widetilde{\eta}$ with
\begin{align*}
 \widetilde{\eta}(\omega,\nu,\kappa) &= F(\nu+s\kappa)+G(\nu-s\kappa)
+H(\omega+\kappa), \\
 \widetilde{q}(\omega,\nu,\kappa) &= s F(\nu+s\kappa)-sG(\nu-s\kappa),
\end{align*}
where $F,G,H$ are (arbitrary) convex functions.\\

The initial conditions  \eqref{datoincial} becomes 
\begin{equation} \label{datoinicialLin}
\begin{cases}
 (\omega(0,y),\nu(0,y),\kappa(0,y))=(\omega_0(y),\nu_0(y),\kappa_0(y)), \quad y \in \re,\\
 \omega_0(y) \ge \underline{\omega} > 0.
\end{cases}
\end{equation}
Due to the fact that system \eqref{sistemaLin} is linear, 
the explicit solution of the corresponding Cauchy problem \eqref{sistemaLin}--\eqref{datoinicialLin} is
\begin{equation}
 \begin{aligned}
 \omega(t,y) &= \omega_0(y)+\kappa_0(y)-\kappa(t,y),\\
 \nu(t,y) &= \frac12 (\nu_0(y+st)+\nu_0(y-st))-
 \frac{s}{2} (\kappa_0(y+st)-\kappa_0(y-st)),\\
 \kappa(t,y) &= \frac12 (\kappa_0(y+st)+\kappa_0(y-st))-
 \frac{1}{2s} (\nu_0(y+st)-\nu_0(y-st)).
\end{aligned}
\end{equation}
Moreover, by condition H1 we obtain that
\begin{align*}
 c_1 \le \nu(t,y)-s\kappa(t,y) \le c_2, \quad c_3 \le \nu(t,y)+s\kappa(t,y) \le c_4,\\
 \omega(t,y)+\kappa(t,y) > c_5,
\end{align*}
and since $\rho_0(x) \ge \underline{\rho}=const.>0$ by \eqref{datoincial}, we have that
$\omega(t,y) \ge \underline{\omega} >0$, ensuring that the function $y \mapsto X(t,y)$ is invertible
and bi-Lipschitzian from $\re$ to $\re$ for all $t \ge 0$, and by \cite{Wagner, Peng2}, 
we also have uniqueness of the entropy solution of \eqref{sistema}--\eqref{datoincial} 
if and only if we have uniqueness of the entropy solution of \eqref{sistemaLin}--\eqref{datoinicialLin}.\\
Therefore, we consider $X_0=Y_0^{-1}$. Then, the unique function $x=X(t,y)$ that satisfy $X(0,y)=X_0(y)$ is given by
\begin{equation} \label{EcuacionA1}
\begin{aligned}
 X(t,y)= & \frac{1}{2s} \int_{y-st}^{y+st} u_0(X_0(\xi)) \, d\xi + \int_0^{y} \left( v_0(X_0(\xi))+\frac{1}{\rho_0(X_0(\xi))} \right) \, d\xi -\\
 & - \frac12 \int_0^{y+st} v_0(X_0(\xi)) \, d\xi - \frac12 \int_0^{y-st} v_0(X_0(\xi)) \, d\xi.
\end{aligned}
\end{equation}
From the above, we obtain the following Theorem.
\begin{teor} \label{Teor1.1}
Assume that $\rho_0,u_0,v_0 \in L^\infty(\re)$ with $\rho_0(x) \ge \underline{\rho}>0$, the conditions H1, H2 hold and
\begin{equation}
\inf_{x \in \re} \left( u_0(x)+\frac{s}{\rho_0(x)} \right) > \sup_{x \in \re} \left( u_0(x)-\frac{s}{\rho_0(x)} \right).
\end{equation}
Then, the Cauchy problem  \eqref{sistema}--\eqref{datoincial}
has an unique global solution $(\rho,u,v)\in L^\infty (\re^+ \times \re)$ that satisfy the entropy condition  \eqref{condEntropia} 
for all pair  $(\eta,q)$ defined in  \eqref{entropia}--\eqref{flujoentropia}. Moreover, this solution is given by
\begin{align*}
 \rho(t,x)  =& \frac{\rho_0(X_0(Y(t,x)))}{1+\rho_0(X_0(Y(t,x))) \left[ v_0(X_0(Y(t,x))) -v(t,x) \right] }, \\
 u(t,x) = & \Gamma^+(t,x) -s \Upsilon^-(t,x) \text{ and } \\
 v(t,x)=& \Upsilon^+(t,x)-\frac{1}{s}\Gamma^-(t,x),
\end{align*}
where
\begin{align*}
 \Gamma^{\pm}(t,x) &= \frac{1}{2} \left[ u_0(X_0(Y(t,x)+st)) \pm u_0(X_0(Y(t,x)-st)) \right]
\end{align*}
and
\begin{align*}
 \Upsilon^{\pm}(t,x) &= \frac{1}{2} \left[ v_0(X_0(Y(t,x)+st)) \pm v_0(X_0(Y(t,x)-st)) \right].
\end{align*}
\end{teor}
Now we apply our result to particular example, which behaves as the advection equation.\\

\emph{Example.} We consider the initial data $u_0(x)=\overline{u}$ and $v_0(x)=\overline{v}$  as being constant functions, 
and $\rho_0(x) \ge 0$ as a bounded function.\\ By the E-L transformation,
$$ Y_0(x)=\int_0^x \rho_0(\xi) \, d\xi , \quad \text{and consider } X_0=Y_0^{-1}. $$
Then,
$$ X(t,y)=X_0(y)+\overline{u}t, \quad \text{and} \quad Y(t,x)=Y_0(x-\overline{u}t). $$
In this way the solution of the Cauchy problem is given by
$$ \rho(t,x)=\rho_0(x-\overline{u}t), \quad u(t,x)=\overline{u}, \quad v(t,x)=\overline{v}. $$
%


\section{Riemann problem} \label{referencia}
In this section we study the solution for the Riemann problem associated with the sBB system,  
in which the left and right constant states $(\rho_l,u_l,v_l)$ and $(\rho_r,u_r,v_r)$, respectively, satisfy the conditions
H1-H2 and $\lambda_1(\rho_l,u_l,v_l) < \lambda_3(\rho_r,u_r,v_r)$.\\

Consider the Riemann problem  of the system \eqref{sistema} with initial data
\begin{equation} \label{datoIniRiemann}
 (\rho,u,v)(0,x)=\begin{cases}
                  (\rho_r,u_r,v_r), & \text{if } x >0,\\
                  (\rho_l,u_l,v_l), & \text{if } x <0,
                 \end{cases}
\end{equation}
where $(\rho_0,u_0,v_0)(x)=(\rho,u,v)(0,x)$ satisfies the conditions H1 and H2.\\
First, observe that system \eqref{sistema} is equivalent to
\begin{equation} \label{sistema1}
 \begin{cases}
  \rho_t + m_x=0, \\
  m_t + (\frac{m^2}{\rho}+s^2 \frac{n}{\rho})_x=0, \\
  n_t + (\frac{mn}{\rho} + \frac{m}{\rho})_x=0,
 \end{cases}
\end{equation}
with $s=const. >0$,
where $m=\rho u$, $n=\rho v$, and the initial data \eqref{datoIniRiemann} is given by
\begin{equation} \label{datoIniRiemannNew}
 (\rho,m,n)(0,x)=\begin{cases}
                  (\rho_r,m_r,n_r), & \text{if } x >0,\\
                  (\rho_l,m_l,n_l), & \text{if } x <0,
                 \end{cases}
\end{equation}
with $m_r=\rho_r u_r$, $n_r=\rho_r v_r$, $m_l=\rho_l u_l$ and $n_l=\rho_l v_l$.\\
The eigenvalues of the system \eqref{sistema}, in the variables $\rho,m,n$, are given by
\begin{equation}
 \lambda_1=\frac{m-s}{\rho}, \quad \quad \lambda_2= \frac{m}{\rho}, \quad \quad \lambda_3= \frac{m+s}{\rho},
\end{equation}
the right eigenvectors become
\begin{align}
 r_1(\rho,m,n) &= \frac{(\rho,m-s,n+1)}{\sqrt{\rho^2+(m-s)^2+(n+1)^2}}, \label{vectorR1} \\
 r_2(\rho,m,n) &= \frac{(\rho,m,n)}{\sqrt{\rho^2+m^2+n^2}} \text{ and } \label{vectorR2} \\
 r_3(\rho,m,n) &= \frac{(\rho,m+s,n+1)}{\sqrt{\rho^2+(m+s)^2+(n+1)^2}}. \label{vectorR3}
\end{align}

From  \eqref{vectorR1} we get that the 1-rarefaction curve can be found as,
\begin{align*}
 \frac{d\rho}{dt}&=\frac{\rho}{\sqrt{\rho^2+(m-s)^2+(n+1)^2}}, \quad \quad \rho(0)=\rho_0, \\
 \frac{dm}{dt}&=\frac{m-s}{\sqrt{\rho^2+(m-s)^2+(n+1)^2}}, \quad \quad m(0)=m_0 \text{ and } \\
 \frac{dn}{dt}&=\frac{n+1}{\sqrt{\rho^2+(m-s)^2+(n+1)^2}}, \quad \quad n(0)=n_0.
\end{align*}
That is,
\begin{equation*}
 \sqrt{\rho^2+(m-s)^2+(n+1)^2}-\sqrt{\rho_0^2+(m_0-s)^2+(n_0+1)^2}=t.
\end{equation*}
Therefore, the integral curves of the vector field $r_1$ are given by straight lines in the direction of the vector $r_1(\rho_0,m_0,n_0)$
and goes trough the point $(0,s,-1)$, that is,
\begin{equation}
 \begin{cases}
   \rho = \rho_0 + \frac{\rho_0}{\sqrt{\rho_0^2+(m_0-s)^2+(n_0+1)^2}}t, \\
   m=m_0 + \frac{m_0-s}{\sqrt{\rho_0^2+(m_0-s)^2+(n_0+1)^2}}t, \\
   m=m_0 + \frac{n_0+1}{\sqrt{\rho_0^2+(m_0-s)^2+(n_0+1)^2}}t,
 \end{cases}  \quad \quad t \in \re.
\end{equation}

Analogously, from \eqref{vectorR2} we can analyze the 2-rarefaction curve. In this case the integral curves corresponding to the vector field $r_2$ are given by straight lines 
going through the origin in the direction of the vector $r_2(\rho_0,m_0,n_0)$. 
Also, from \eqref{vectorR3} we can see that for the 3-rarefaction curve, the integral curves of the vector field $r_3$,
are give by straight lines trough  $(0,-s,-1)$ that are parallel to $r_3(\rho_0,m_0,n_0)$.\\
Thus, the i-rarefaction curve $R_i(\sigma)(\rho,m,n)$ satisfy
\begin{equation} \label{conditionRarefaction}
 R_i(\sigma)(\rho,m,n)=(\rho,m,n)+\sigma r_i(\rho,m,n) \text{ for } i=1,2,3.
\end{equation}

This also may be deduced from self-similar solution
\begin{equation}
 (\rho,u,v)(t,x)= (\rho,u,v)(\xi), \quad \xi=\frac{x}{t},
\end{equation}
for which system \eqref{sistema} becomes
\begin{equation}\label{SistemaSelfSol}
 \begin{cases}
  &-\xi \rho_\xi+(\rho u)_\xi=0,\\
  &-\xi (\rho u)_\xi+(\rho u^2+s^2v)_\xi=0,\\
  &-\xi (\rho v)_\xi+(\rho uv+u)_\xi=0,\\
 \end{cases}
\end{equation}
and initial data \eqref{datoIniRiemann} changes to the boundary condition
\begin{equation}
 \begin{aligned}
  (\rho,u,v)(-\infty)=(\rho_l,u_l,v_l) \text{ and } 
  (\rho,u,v)(+\infty)=(\rho_r,u_r,v_r).
 \end{aligned}
\end{equation}
This is a two-point boundary value problem of first-order ordinary differential equations
with the boundary values in the infinity.\\
For smooth solution, \eqref{SistemaSelfSol} is reduced to
\begin{equation}
 \begin{pmatrix}
  u-\xi & \rho & 0 \\ 0 & \rho(u-\xi) & 0 \\ 0 & 1 & \rho(u-\xi)
 \end{pmatrix}
 \begin{pmatrix}
  \rho \\ u \\ v
 \end{pmatrix}_\xi =0.
\end{equation}
It provides either the general solutions (constant states)
\begin{equation}
 (\rho,u,v)=constant \quad \quad (\rho>0),
\end{equation}
or singular solutions
\begin{equation} \label{SolSing1}
 \begin{aligned}
  \xi=\lambda_1=u-\frac{s}{\rho}, \quad \text{d}\left( u-\frac{s}{\rho} \right)=0 \text{ and } \text{d}\left( v+\frac{1}{\rho} \right)=0,\\
  \xi=\lambda_2=u, \quad \text{d}u=0 \text{ and } \text{d}v=0, \\
  \xi=\lambda_3=u+\frac{s}{\rho}, \quad \text{d}\left( u+\frac{s}{\rho} \right)=0 \text{ and } \text{d}\left( v+\frac{1}{\rho} \right)=0.
 \end{aligned}
\end{equation}
Integrating \eqref{SolSing1} from $(\rho_l,u_l,v_l)$ to $(\rho,u,v)$, one can get that
\begin{equation} \label{SolSing2}
 \begin{aligned}
  \xi=\lambda_1=u-\frac{s}{\rho}, \quad u-\frac{s}{\rho} =u_l-\frac{s}{\rho_l}  \text{ and } v+\frac{1}{\rho} =v_l+\frac{1}{\rho_l},\\
  \xi=\lambda_2=u, \quad u=u_l \text{ and } v=v_l, \\
  \xi=\lambda_3=u+\frac{s}{\rho}, \quad u+\frac{s}{\rho} =u_l+\frac{s}{\rho_l} \text{ and } v+\frac{1}{\rho} =v_l+\frac{1}{\rho_l}.
 \end{aligned}
\end{equation}
Oberve that \eqref{SolSing2} in the variables $\rho,m,n$ is equivalent to \eqref{conditionRarefaction}.\\
For a bounded discontinuity at $\xi=\omega$, the Rankine-Hugoniot conditions hold. That is,
\begin{equation} \label{RHcondition}
 \begin{cases}
  -\omega [\rho]+[\rho u]=0,\\
  -\omega [\rho u]+[\rho u^2+s^2v]=0,\\
  -\omega [\rho v]+[\rho uv+u]=0,
 \end{cases}
\end{equation}
where $[q] = q_l-q$ is the jump of $q$ across the discontinuous line and 
$\omega$ is the velocity of the discontinuity.
From \eqref{RHcondition}, we have
\begin{equation} \label{SolRH1}
 \begin{aligned}
  \omega=u-\frac{s}{\rho}, \quad u-\frac{s}{\rho} =u_l-\frac{s}{\rho_l}  \text{ and } v+\frac{1}{\rho} =v_l+\frac{1}{\rho_l},\\
  \omega=u, \quad u=u_l \text{ and } v=v_l, \\
  \omega=u+\frac{s}{\rho}, \quad u+\frac{s}{\rho} =u_l+\frac{s}{\rho_l} \text{ and } v+\frac{1}{\rho} =v_l+\frac{1}{\rho_l}.
 \end{aligned}
\end{equation}
From \eqref{SolSing2} and \eqref{SolRH1}, we conclude that the rarefaction waves and the shock waves are coincident, which correspond to contact discontinuities.
Namely, for a given left state $(\rho_l,u_l,v_l)$, the contact discontinuity curves, which are the sets of states that can be connected on the right by a 1-contact discontinuity $J_1$, a 2-contact discontinuity $J_2$ or a 3-contact discontinuity $J_3$, are as follows:
\begin{equation}\label{discontinuities}
 \begin{aligned}
  &J_1: (\rho,u,v):=(\rho, u_l-s/\rho_l+s/\rho,v_l+1/\rho_l-1/\rho),\\
  &J_2: (\rho,u,v):=(\rho, u_l,v_l),\\
  &J_3: (\rho,u,v):=(\rho, u_l+s/\rho_l-s/\rho,v_l+1/\rho_l-1/\rho),\quad \rho >0.
 \end{aligned}
\end{equation}
In the space $(\rho > 0, u\in \re,u\in \re)$, through the point $(\rho_l,u_l,v_l)$, we draw curves \eqref{discontinuities} which are denoted by $J_1$, $J_2$ and $J_3$ respectively. So, $J_1$ has asymptotes $\rho=0$  and $(\rho, u_l-s/\rho_l,v_l+1/\rho_l)$ for $\rho \ge 0$, and $J_3$ has asymptotes $\rho=0$ and $(\rho,u_l+s/\rho_l,v_l+1/\rho_l)$.

In order to solve the Riemann problem \eqref{sistema}--\eqref{datoIniRiemann},
we consider constant left and right states $U_l=(\rho_l,u_l,v_l)$ and $U_r=(\rho_r,u_r,v_r)$, respectively, such that
the conditions H1-H2 are satisfied and $\lambda_1(U_l)<\lambda_3(U_r)$. Then there exists intermediate states,
 $U_*=(\rho_*,u_*,v_*)$
and $U_{**}=(\rho_{**},m_{**},n_{**})$ such that $U_*=R_1(\sigma_1)(U_l)$, $U_{**}=R_2(\sigma_2)(U_*)$ and  $U_r=R_3(\sigma_3)(U_{**})$,
for some $\sigma_1$, $\sigma_2$ y $\sigma_3$.\\
Furthermore, because of  \eqref{conditionRarefaction}, the states $U^*$ y $U^{**}$ should satisfy
\begin{subequations} \label{eq34}
\begin{align}
  u_*&=\left( u_l-\frac{s}{\rho_l} \right)+\frac{s}{\rho_*}, \label{for34a} \\
  v_*&=\left( v_l+\frac{1}{\rho_l} \right)-\frac{1}{\rho_*}, \label{for34b} \\
  u_*&=u_{**}, \label{for34c}\\
  v_*&=v_{**}, \label{for34d}\\
  u_{**}&=\left( u_r+\frac{s}{\rho_r} \right)-\frac{s}{\rho_{**}} \text{ and } \label{for34e} \\
  v_{**}&=\left( v_r+\frac{1}{\rho_r} \right)-\frac{1}{\rho_{**}}. \label{for34f}
\end{align}
\end{subequations}
Now, we need to ensure that $\rho_*$ and $\rho_{**}$ be positive. Because if $\rho_*$ or $\rho_{**}$ are negative 
the problem has mathematical interest, but not physical.
From equations  \eqref{for34d}, \eqref{for34b} and \eqref{for34f} we have
\begin{equation} \label{eq35}
 \frac{1}{\rho_{**}}-\frac{1}{\rho_*}= \left( v_r+\frac{1}{\rho_r} \right) - \left( v_l+\frac{1}{\rho_l} \right).
\end{equation}
On the other hand, from \eqref{for34c}, \eqref{for34a} and \eqref{for34e}, we have
\begin{equation} \label{eq36}
 \frac{1}{\rho_{**}}+\frac{1}{\rho_*}=\frac{1}{s} \left\{ \left( u_r+\frac{s}{\rho_r} \right) - \left( u_l-\frac{s}{\rho_l} \right) \right\}.
\end{equation}
From  \eqref{eq35} y \eqref{eq36} we conclude that
\begin{align}
 \frac{1}{\rho_*}= \frac{1}{2s} \left( \lambda_3(U_r) - \lambda_1(U_l) \right) -\frac12 \left( R_2(U_r)- R_2(U_l) \right), \label{eq37} \\
 \frac{1}{\rho^{**}}= \frac{1}{2s} \left( \lambda_3(U_r)-\lambda_1(U_l) \right) +\frac12 \left( R_2(U_r)-R_2(U_l) \right). \label{eq38}
\end{align}
Thus, \eqref{eq34}, \eqref{eq37} and \eqref{eq38} gives
\begin{align}
 u_*=\frac12 \{ (u_l+sv_l)+(u_r-sv_r) \}=u_{**}, \\
 v_*=\frac{1}{2s} \{ (u_l+sv_l) - (u_r-sv_r) \}=v_{**}.
\end{align}
 Observe that by conditions H1 and H2, we have that $U_*$ and $U_{**}$ also satisfies.
This guarantees that $\rho_*$ and $\rho_{**}$ are positive.

Note that if $\lambda_1(U^-)<\lambda_3(U^+)$, then,
\begin{equation*}
 \left| R_2(U^+)-R_2(U^-) \right| < \frac{1}{s} (\lambda_3(U^+)-\lambda_1(U^-)).
\end{equation*}

\begin{nota}
Note that if  $\lambda_1(U_l)=\lambda_3(U_r)$, then
\begin{equation} \label{eq40}
 s=\frac{\frac{m_l}{\rho_l}-\frac{m_r}{\rho_r}}{\frac{1}{\rho_r}+\frac{1}{\rho_l}},
\end{equation}
If we assume that  $\rho_r=\rho_l$ and $n_r=n_l$, then, \eqref{eq40} reduces to $s=\frac{m_l-m_r}{2}$. 
In this case the Riemann problem does not have a solution since 
the lines $R_1(\sigma)(U_l)$ and $R_3(\sigma)(U_{**})$ are parallel.
\end{nota}

Thus, by conditions H1 and H2, for $s>0$ and if $\lambda_3(U_r)>\lambda_1(U_l)$,
the solution of the Riemann problem is given by\begin{equation}\label{solucionRiemann}
 (\rho,u,v)(t,x)=\begin{cases}
                            (\rho_l,u_l,v_l), & \quad \text{if } \frac{x}{t} < \lambda_1(U_l),\\
                            (\rho_*,u_*,v_*), & \quad \text{if } \lambda_1(U_l) \le \frac{x}{t} < \lambda_2(U_*),\\
                            (\rho_{**},u_{**},v_{**}), & \quad \text{if } \lambda_2(U_{**}) \le \frac{x}{t} < \lambda_3(U_r),\\
                            (\rho_r,u_r,v_r), & \quad \text{if } \frac{x}{t} \ge \lambda_3(U_r).
                           \end{cases}
\end{equation}
Additionally, as usual, since the system is linearly degenerate,  $\lambda_1(U_l)=\lambda_1(U_*)$, $\lambda_2(U_*)=\lambda_2(U_{**})$ and $\lambda_3(U_{**})=\lambda_3(U_r)$.

\begin{lema}
 Given left and right constant states $(\rho_l,u_l,v_l)$ and $(\rho_r,u_r,v_r)$, respectively, such that they
 satisfy conditions H1,H2 and
 $\lambda_1(\rho_l,u_l,v_l)<\lambda_3(\rho_r,u_r,v_r)$. Then, the following relation is satisfied
 \begin{equation}
  | R_2(\rho_r,u_r,v_r)-R_2(\rho_l,u_l,v_l) | < \frac{1}{s} (\lambda_3(\rho_r,u_r,v_r)-\lambda_1(\rho_l,u_l,v_l)).
 \end{equation}
\end{lema}

The results of this section can be summarized in the following Theorem.

\begin{teor}
 Given left and right constant states $(\rho_l,u_l,v_l)$ and $(\rho_r,u_r,v_r)$, respectively, such that they satisfy 
 conditions H1,H2 and $\lambda_1(\rho_l,u_l,v_l)<\lambda_3(\rho_r,u_r,v_r)$.
Then, there is a unique global solution to the Riemann problem \eqref{sistema}--\eqref{datoIniRiemann}.
Moreover, this solution is given by
\begin{equation}
 (\rho,u,v)(t,x)=\begin{cases}
                  (\rho_r,u_r,v_r), & \text{if } x\ge \lambda_3(\rho_r,u_r,v_r), \\
                  (\rho_{**},u_{**},v_{**}), & \text{if }  \lambda_2(\rho_{**},u_{**},v_{**}) \le x < \lambda_3(\rho_r,u_r,v_r), \\
                  (\rho_*,u_*,v_*), & \text{if } \lambda_1(\rho_l,u_l,v_l) \le x < \lambda_2(\rho_{**},u_{**},v_{**}), \\
                  (\rho_l,u_l,v_l), & \text{if } x < \lambda_1(\rho_l,u_l,v_l),
                 \end{cases}
\end{equation}
where
\begin{align}
 \frac{1}{\rho_*}&=\frac{1}{2s}(u_r-u_l)-\frac12 (v_r-v_l)+\frac{1}{\rho_l}, \\
 \frac{1}{\rho_{**}}&=\frac{1}{2s}(u_r-u_l)+\frac12 (v_r-v_l)+\frac{1}{\rho_r}, \\
 u_*&=\frac12 \{ (u_l+sv_l)+(u_r-sv_r) \}=u_{**} \text{ and } \\
 v_*&=\frac{1}{2s} \{ (u_l+sv_l) - (u_r-sv_r) \}=v_{**}.
\end{align}
\end{teor}
\nocite{bressan1}

\section{Delta shock solution}

In this section, we discuss the solution for the Riemann problem associated with the sBB system, in which the left and right constant
states $(\rho_l,u_l,v_l)$ and $(\rho_r,u_r,v_r)$, respectively, satisfy the conditions H1 and H2, but unlike previous section they
satisfy $\lambda_1(\rho_l,u_l,v_l) \ge \lambda_3(\rho_r,u_r,v_r)$.\\

Denote by $BM(\re)$ the space of bounded Borel measures on $\re$, and then the definition of a measure solution of sBB system in $BM(\re)$ can be given as follows.
\renewcommand{\labelenumi}{\alph{enumi}$)$ }
\begin{defin} \label{measureSol}
 A triple $(\rho,u,v)$ constitutes a \emph{measure solution} to the sBB system, if it holds that
 \begin{enumerate}
  \item $\rho \in L^\infty((0,\infty),BM(\re)) \cap C((0,\infty),H^{-s}(\re))$,
  \item $u \in L^\infty((0,\infty),L^\infty(\re)) \cap C((0,\infty),H^{-s}(\re))$,
  \item $v \in L_{loc}^\infty((0,\infty),L_{loc}^\infty(\re)) \cap C((0,\infty),H^{-s}(\re))$, \ \ $s>0$,
  \item $u$ and $v$ are measurable with respect to $\rho$ at almost for all $t\in(0,\infty)$,
 \end{enumerate}
and
\begin{equation} \label{MS1A}
\begin{cases}
 I_1=\int_0^\infty \int_{\re} (\phi_t+u\phi_x) \, d\rho dt =0,\\
 I_2=\int_0^\infty \int_{\re} u(\phi_t+u\phi_x) \, d\rho dt + \int_0^\infty \int_{\re} s^2 v \phi_x \, dx dt =0, \\
 I_3=\int_0^\infty \int_{\re} v(\phi_t+u\phi_x) \, d\rho dt + \int_0^\infty \int_{\re} u \phi_x \, dx dt =0,
\end{cases}
\end{equation}
for all test function $\phi \in C_0^\infty(\re^+ \times \re)$.
\end{defin}

\begin{defin}
 A two-dimensional weighted delta function $w(s)\delta_L$ supported on a smooth curve $L$ parameterized as $t=t(s)$, $x=x(s)$ $(c \le s \le d)$ is defined by
 \begin{equation}
  \langle w(s)\delta_L, \phi(t,x) \rangle = \int_c^d w(s)\phi(t(s),x(s)) \, ds
 \end{equation}
for all $\phi \in C_0^\infty(\re^2)$.
\end{defin}

\begin{defin}
 A triple distribution $(\rho,u,v)$ is called a \emph{delta shock wave} if it is represented in the form
 \begin{equation} \label{deltashock}
  (\rho,u,v)(t,x)=\begin{cases}
                   (\rho_l,u_l,v_l)(t,x), & x<x(t),\\
                   (w(t)\delta(x-x(t)),u_\delta(t),g(t)), & x=x(t),\\
                   (\rho_r,u_r,v_r)(t,x), & x>x(t),
                  \end{cases}
 \end{equation}
and satisfies Definition \ref{measureSol}, where $(\rho_l,u_l,v_l)(t,x)$ and $(\rho_r,u_r,v_r)(t,x)$ are piecewise smooth bounded solutions of 
the sBB system \eqref{sistema}.
\end{defin}

We set $\frac{dx}{dt}=u_\delta(t)$ since the concentration in $\rho$ need to travel at the speed of discontinuity. Hence, we say that a delta shock wave \eqref{deltashock} is a measure solution to the sBB system \eqref{sistema} if and only if the
following relation holds,
\begin{equation} \label{RHgener}
 \begin{cases}
  \frac{dx(t)}{dt}=u_\delta(t),\\
  \frac{dw(t)}{dt}=-[\rho]u_\delta(t)+[\rho u],\\
  \frac{dw(t)u_\delta(t)}{dt}=-[\rho u]u_\delta(t)+[\rho u^2+s^2v],\\
  \frac{dw(t)g(t)}{dt}=-[\rho v]u_\delta(t)+[\rho uv+u].
 \end{cases}
\end{equation}
In fact, for any test function $\phi \in C_0^\infty (\re^+ \times \re)$,
from \eqref{MS1A}, we obtain
\begin{equation*}
 \begin{aligned}
  I_1 &=\int_0^\infty \int_{\re} (\phi_t+u\phi_x) \, d\rho dt =\int_0^\infty \left\{ -u_\delta(t)[\rho]+[\rho u]-\frac{dw(t)}{dt} \right\} \phi \, dt,\\
 I_2 &=\int_0^\infty \int_{\re} u(\phi_t+u\phi_x) \, d\rho dt + \int_0^\infty \int_{\re} s^2 v \phi_x \, dx dt \\ &= \int_0^\infty \left\{ -u_\delta(t)[\rho u]+[\rho u^2+s^2v]-\frac{dw(t)u_\delta(t)}{dt} \right\} \, dt, \quad \text{ and } \\
 I_3 &= \int_0^\infty \int_{\re} v(\phi_t+u\phi_x) \, d\rho dt + \int_0^\infty \int_{\re} u \phi_x \, dx dt \\ 
 &=\int_0^\infty \left\{ -u_\delta(t)[\rho v]+[\rho uv+u]-\frac{dw(t)g(t)}{dt} \right\} \phi \, dt.
 \end{aligned}
\end{equation*}
Relations \eqref{RHgener} is called the generalized Rankine-Hugoniot relation. It  reflects the exact relationship among the
limit states on two sides of the discontinuity, the weight,
propagation speed and the location of the discontinuity.
In addition, to guarantee uniqueness, the delta shock wave should satisfy 
the admissibility (entropy) condition
\begin{equation}
 \lambda_3(\rho_r,u_r,v_r) \le u_\delta(t) \le \lambda_1(\rho_l,u_l,v_l).
\end{equation}

Now, the generalized Rankine-Hugoniot relation is applied to the Riemann problem \eqref{sistema}--\eqref{datoIniRiemann} with left and right
 constant states $U_-=(\rho_-,u_-,v_-)$ and
$U_+=(\rho_+,u_+,v_+)$, respectively, satisfying the conditions H1 and H2, the fact $\lambda_3(\rho_+,u_+,v_+) \le \lambda_1(\rho_-,u_-,v_-)$ and
\begin{equation} \label{CondicionDeltaSol}
\begin{aligned}
 \frac12 & (\lambda_1(U_-)-\lambda_3(U_+))^2 \ge \\ & \max \left\{ - \frac{s^2}{\rho_+}(R_2(U_+)-R_1(U_-)), \frac{s^2}{\rho_-}(R_2(U_+)-R_1(U_-)) \right\}.
\end{aligned}
\end{equation}

Thereby, the Riemann problem is reduced to solving \eqref{RHgener} with initial data
\begin{equation} \label{datoRHgener}
 t=0, \quad x(0)=0, w(0)=0, g(0)=0,
\end{equation}
under entropy condition
\begin{equation} \label{condEntroGen}
 u_++\frac{s}{\rho_+} \le u_\delta(t) \le u_--\frac{s}{\rho_-}.
\end{equation}
From \eqref{RHgener} and \eqref{datoRHgener}, it follows that
\begin{equation}\label{DFREW}
 \begin{aligned}
  w(t) &= -[\rho]x(t)+[\rho u]t, \\
  w(t)u_\delta(t) &= -[\rho u]x(t)+[\rho u^2+s^2v]t, \text{  and } \\
  w(t)g(t) &= -[\rho v]x(t)+[\rho uv+u]t.
 \end{aligned}
\end{equation}
Multiplying the first equation in \eqref{DFREW} by $u_\delta(t)$ and then subtracting it from the second one, we obtain that
\begin{equation}
 [\rho]x(t)u_\delta(t)-[\rho u]u_\delta(t)t-[\rho u]x(t)+[\rho u^2+s^2v]t=0,
\end{equation}
that is,
\begin{equation}
 \frac{d}{dt} \left( \frac{[\rho]}{2}x^2(t)-[\rho u]x(t)t+\frac{[\rho u^2+s^2v]}{2}t^2 \right) =0,
\end{equation}
which is equivalent to
\begin{equation} \label{cuadratica1}
 [\rho]x^2(t)-2[\rho u]x(t)t+[\rho u^2+s^2v]t^2 =0.
\end{equation}
From \eqref{cuadratica1}, one can find $u_\delta(t) := u_\delta$ is a constant and $x(t) = u_\delta t$. Then, \eqref{cuadratica1} can
be rewritten
\begin{equation} \label{cuadratica2}
 [\rho]u_\delta^2-2[\rho u]u_\delta+[\rho u^2+s^2v]=0.
\end{equation}

When $[\rho]=\rho_--\rho_+=0$,the situation is very simple and one can easily calculate the solution
\begin{equation} \label{SolDeltaLin}
 \begin{cases}
  u_\delta=\frac{u_-+u_+}{2}+s^2\frac{[v]}{2\rho_-[u]}, \\
  x(t)= u_\delta t, \\
  w(t)=\rho_-(u_--u_+)t, \\
  g(t)= \frac{[\rho u v + u]-u_\delta}{[\rho u]},
 \end{cases}
\end{equation}
which obviously satisfies the entropy condition \eqref{condEntroGen}, 	
since by condition \eqref{CondicionDeltaSol},
\begin{align*}
 s^2 \frac{[v]}{\rho_-} \le \frac12 (\lambda_1(U_-)-\lambda_3(U_+))^2 < \frac12 [u] (\lambda_1(U_-)-\lambda_3(U_+))
\end{align*}
and
\begin{align*}
 u_\delta-\left(u_--\frac{s}{\rho_-} \right) &= \frac{u_-+u_+}{2}+s^2\frac{[v]}{2\rho_-[u]} -\left(u_--\frac{s}{\rho_-} \right) \\
 & =\frac12 \left(  \left( u_++\frac{s}{\rho_-} \right) - \left( u_--\frac{s}{\rho_-} \right) + s^2 \frac{[v]}{\rho_-[u]} \right) \le 0.
\end{align*}
Similarly we can deduce that
\begin{align*}
 u_\delta-\left(u_++\frac{s}{\rho_-} \right) &= \frac{u_-+u_+}{2}+s^2\frac{[v]}{2\rho_-[u]} -\left(u_++\frac{s}{\rho_-} \right) \\
 & =\frac12 \left(  \left( u_--\frac{s}{\rho_-} \right) - \left( u_++\frac{s}{\rho_-} \right) + s^2 \frac{[v]}{\rho_-[u]} \right) \ge 0,
\end{align*}
because
\begin{align*}
 -s^2 \frac{[v]}{\rho_-} \le \frac12 (\lambda_1(U_-)-\lambda_3(U_+))^2 < \frac12 [u] (\lambda_1(U_-)-\lambda_3(U_+)).
\end{align*}
When $[\rho]=\rho_--\rho_+ \neq 0$, the discriminant of the quadratic equation \eqref{cuadratica2} is
\begin{equation}
 \Delta = 4[\rho u]^2-4[\rho][\rho u^2+s^2v]=\rho_-\rho_+[u]^2-s^2[\rho][v]>0
\end{equation}
and then we can find
\begin{equation} \label{SolDelta1}
 \begin{cases}
  u_\delta=\frac{[\rho u]-\sqrt{[\rho u]^2-[\rho][\rho u^2+s^2 v]}}{[\rho]},\\
  x(t)=\frac{[\rho u]-\sqrt{[\rho u]^2-[\rho][\rho u^2+s^2 v]}}{[\rho]} t,\\
  w(t)= \sqrt{[\rho u]^2-[\rho][\rho u^2+s^2 v]} t,\\
  g(t)= \frac{-[\rho u][\rho v]+[\rho v]\sqrt{[\rho u]^2-[\rho][\rho u^2+s^2 v]}+[\rho][\rho u v +u]}{[\rho]\sqrt{[\rho u]^2-[\rho][\rho u^2+s^2 v]}}t,
 \end{cases}
\end{equation}
or, 
\begin{equation} \label{SolDelta2}
 \begin{cases}
  u_\delta=\frac{[\rho u]+\sqrt{[\rho u]^2-[\rho][\rho u^2+s^2 v]}}{[\rho]},\\
  x(t)=\frac{[\rho u]+\sqrt{[\rho u]^2-[\rho][\rho u^2+s^2 v]}}{[\rho]} t,\\
  w(t)= -\sqrt{[\rho u]^2-[\rho][\rho u^2+s^2 v]} t,\\
  g(t)= \frac{-[\rho u][\rho v]-[\rho v]\sqrt{[\rho u]^2-[\rho][\rho u^2+s^2 v]}+[\rho][\rho u v +u]}{[\rho]\sqrt{[\rho u]^2-[\rho][\rho u^2+s^2 v]}}t.
 \end{cases}
\end{equation}

Next, with the help of the entropy condition \eqref{condEntroGen}, we will choose the admissible solution from \eqref{SolDelta1} and \eqref{SolDelta2}. 
Observe that by the entropy condition and since the system is strictly hyperbolic, we have that
$$ u_+-\frac{s}{\rho_+}<u_+<u_++\frac{s}{\rho_+} \le u_--\frac{s}{\rho_-} <u_-<u_-+\frac{s}{\rho_-}. $$
Observe that,
\begin{align*}
 -[\rho]\lambda_1(\rho_-,u_-,v_-)+[\rho u]=\rho_+ \left( \left(u_--\frac{s}{\rho_-} \right) - \left(u_+-\frac{s}{\rho_+} \right) \right) >0,\\
 -[\rho]\lambda_3(\rho_+,u_+,v_+)+[\rho u]=\rho_- \left( \left(u_-+\frac{s}{\rho_-} \right) - \left(u_++\frac{s}{\rho_+} \right) \right) >0,
\end{align*}
\begin{align*}
 [\rho](\lambda_1(\rho_-,u_-,v_-))^2-2[\rho u]\lambda_1(\rho_-,u_-,v_-)+[\rho u^2]+s^2[v] = \\
 -\rho_+\left( u_--u_+-\frac{s}{\rho_-} \right)^2+\frac{s^2}{\rho_-}+s^2[v] \le 0,
\end{align*}
\begin{align*}
 [\rho](\lambda_3(\rho_+,u_+,v_+))^2-2[\rho u]\lambda_3(\rho_+,u_+,v_+)+[\rho u^2]+s^2[v] = \\
 \rho_-\left( u_--u_+-\frac{s}{\rho_+} \right)^2-\frac{s^2}{\rho_+}+s^2[v] \ge 0,
\end{align*}
then, for the solution given in \eqref{SolDelta1}, we have
\begin{align*}
 u_\delta- & \lambda_1(\rho_-,u_-,v_-)= \\ & \frac{[\rho](\lambda_1(\rho_-,u_-,v_-))^2-2[\rho u]\lambda_1(\rho_-,u_-,v_-)+[\rho u^2]+s^2[v]}{(-[\rho]\lambda_1(\rho_-,u_-,v_-)+[\rho u])+\sqrt{[\rho u]^2-[\rho][\rho u^2+s^2 v]}} \le 0
\end{align*}
and
\begin{align*}
 u_\delta- & \lambda_3(\rho_+,u_+,v_+)= \\ & \frac{[\rho](\lambda_3(\rho_+,u_+,v_+))^2-2[\rho u]\lambda_3(\rho_+,u_+,v_+)+[\rho u^2]+s^2[v]}{(-[\rho]\lambda_3(\rho_+,u_+,v_+)+[\rho u])+\sqrt{[\rho u]^2-[\rho][\rho u^2+s^2 v]}} \ge 0,
\end{align*}
which imply that the entropy condition \eqref{condEntroGen} is valid.
When $\lambda_1(\rho_-,u_-,v_-)=\lambda_3(\rho_+,u_+,v_+)$, we have trivially that
$\lambda_1(\rho_-,u_-,v_-)=u_\delta=\lambda_3(\rho_+,u_+,v_+)$.

Now, for the solution \eqref{SolDelta2}, when $\rho_-<\rho_+$ we have
\begin{equation*}
 \begin{aligned}
  u_\delta &-\lambda_3(\rho_+,u_+,v_+)=\frac{-[\rho]\lambda_3(U_+)+[\rho u]+\sqrt{[\rho u]^2-[\rho][\rho u^2+s^2v]}}{[\rho]}\\
  &=\frac{\rho_-(\lambda_3(U_-)-\lambda_3(U_+))+\sqrt{[\rho u]^2-[\rho][\rho u^2+s^2v]}}{[\rho]} <0,
 \end{aligned}
\end{equation*}
and when $\rho_->\rho_+$, that
\begin{equation*}
 \begin{aligned}
  u_\delta &-\lambda_1(\rho_-,u_-,v_-)=\frac{-[\rho]\lambda_1(U_-)+[\rho u]+\sqrt{[\rho u]^2-[\rho][\rho u^2+s^2v]}}{[\rho]}\\
  &=\frac{\rho_+(\lambda_1(U_-)-\lambda_1(U_+))+\sqrt{[\rho u]^2-[\rho][\rho u^2+s^2v]}}{[\rho]} >0.
 \end{aligned}
\end{equation*}
showing that the solution \eqref{SolDelta2} does not satisfy the entropy
condition \eqref{condEntroGen}.\\
Thus we have proved the following result.
\begin{teor}
 Given they left and right constant states $(\rho_l,u_l,v_l)$ and $(\rho_r,u_r,v_r)$,
respectively, such that satisfy the conditions H1 and H2, $\lambda_1(\rho_l,u_l,v_l) \ge \lambda_3(\rho_r,u_r,v_r)$
and \eqref{CondicionDeltaSol}, that is,
\begin{equation*}
 \frac12 (\lambda_1(U_l)-\lambda_3(U_r))^2 \ge \max \left\{ -\frac{s^2}{\rho_r}(R_2(U_r)-R_2(U_l)),\frac{s^2}{\rho_l}(R_2(U_r)-R_2(U_l)) \right\}.
\end{equation*}
Then, the Riemann problem \eqref{sistema}--\eqref{datoIniRiemann} admits a unique entropy solution in the sense of measures.
This solution is of the form
\begin{equation}
 (\rho,u,v)(t,x)=\begin{cases}
                  (\rho_l,u_l,v_l), & \text{if } x<u_\delta t, \\
                  (w(t)\delta(x-u_\delta t),u_\delta,g(t)), & \text{if } x=u_\delta t, \\
                  (\rho_r,u_r,v_r), & \text{if } x>u_\delta t,
                 \end{cases}
\end{equation}
where $u_\delta$, $w(t)$ and $g(t)$ are show in \eqref{SolDeltaLin} for $[\rho]=0$ or \eqref{SolDelta1} for $[\rho] \neq 0$.
\end{teor}

We are also interested in studying the case
\begin{equation} \label{eqDel37}
 \frac{1}{s} \left( \lambda_3(U_r) - \lambda_1(U_l) \right) - \left( R_2(U_r)- R_2(U_l) \right)=0,
\end{equation}
or
\begin{equation} \label{eqDel38}
 \frac{1}{s} \left( \lambda_3(U_r)-\lambda_1(U_l) \right) + \left( R_2(U_r)-R_2(U_l) \right)=0, 
\end{equation}
that according to our knowledge is the best case of the delta shock waves (see Equations \eqref{eq37} and \eqref{eq38}).\\

We began by analyzing the case in \eqref{eqDel37}. Assume left and right constant states $(\rho_l,u_l,v_l)$ and $(\rho_r,u_r,v_r)$, 
respectively, such that satisfy the conditions H1 and H2, $\lambda_1(\rho_l,u_l,v_l) \ge \lambda_3(\rho_r,u_r,v_r)$
and \eqref{eqDel37}.\\
It is easy to see that if $\lambda_1(\rho_l,u_l,v_l) = \lambda_3(\rho_r,u_r,v_r)$, then, the inequality \eqref{CondicionDeltaSol} is trivially satisfied.
Suppose that $\lambda_1(\rho_l,u_l,v_l) > \lambda_3(\rho_r,u_r,v_r)$. Then, $\lambda_3(\rho_r,u_r,v_r) < \lambda_3(\rho_l,u_l,v_l)$ , and so
\begin{equation*}
 \lambda_3(\rho_r,u_r,v_r)-\lambda_1(\rho_l,u_l,v_l)< \frac{2s}{\rho_l},
\end{equation*}
meaning that \eqref{CondicionDeltaSol} is satisfied.
The analysis is similar for \eqref{eqDel38}.


\section*{Acknowledgements}
The authors are grateful with Professor Yunguang Lu for suggesting this problem, pointing 
out the simplification of the system proposed by Bouchut and Boyaval and also for 
communicating to us  the preprint \cite{Lu}.


\bibliography{bibliogra}
\bibliographystyle{abbrv}

\end{document}